\documentclass[reqno, 12pt]{amsart}

\usepackage{amssymb}
\usepackage{epsfig}
\usepackage{epsf}
\usepackage{float}
\usepackage{a4wide}

\newcommand{\mc}[1]{{\mathcal #1}}

\newcommand{\bb}[1]{{\mathbb #1}}

\numberwithin{equation}{section}

\renewcommand\thefigure{\thesection.\@arabic\c@figure}
\renewcommand\thetable{\thesection.\@arabic\c@table}

\newtheorem{thm}{Theorem}[section]
\newtheorem{problem}{Problem}[section]
\newtheorem{lemma}[thm]{Lemma}
\newtheorem{corol}[thm]{Corollary}
\newtheorem{propos}[thm]{Proposition}

\def\reff#1{(\ref{#1})}
\usepackage{hyperref}

\begin{document}

\def\E{{\mathbb E}}
\def\P{{\mathbb P}}
\def\R{{\mathbb R}}
\def\Z{{\mathbb Z}}
\def\N{{\mathbb N}}
\def\cF{{\mathcal F}}
\def\cS{{\mathcal S}}
\def\cO{{\mathcal O}}
\def\W{{\cal W}}
\def\G{{\mathcal T}}
\def\T{{\cal T}}
\def\I{{\cal I}}
\def\TT{\bar{{\cal T}}}
\def\II{\bar{{\cal I}}}
\def\C{{\C}}
\def\C{{\cal D}}
\def\n{{\bf n}}
\def\m{{\bf m}}
\def\b{{\bf b}}
\def\Var{{\hbox{Var}}}
\def\Cov{{\hbox{Cov}}}
\def\uR{{\underline R}}
\def\oR{{\overline R}}
\def\urho{{\underline \rho}}
\def\orho{{\overline \rho}}

\def\sqr{\vcenter{
         \hrule height.1mm
         \hbox{\vrule width.1mm height2.2mm\kern2.18mm\vrule width.1mm}
         \hrule height.1mm}}                  % This is a slimmer sqr.
\def\square{\ifmmode\sqr\else{$\sqr$}\fi}
\def\one{{\bf 1}\hskip-.5mm}
\def\limn{\lim_{N\to\infty}}
\def\given{\ \vert \ }
\def\ze{{\zeta}}
\def\be{{\beta}}
\def\la{{\lambda}}
\def\ga{{\gamma}}
\def\a{{\alpha}}
\def\th{{\theta}}
\def\proof{\noindent{\bf Proof. }}
\def\A{{\bf A}}
\def\B{{\bf B}}
\def\C{{\bf C}}
\def\D{{\bf D}}
\def\MM{{\bf m}}
\def\w{\bar{w}}
\def\lnt{{\Lambda^N}}
\def\dlnt{\delta\Lambda^N_t}
\def\lno{\Lambda^N_0}
\def\dlno{\delta\Lambda^N_0}

\title{Poisson trees, succession lines and coalescing random walks}
\date{}

\author{P. A. Ferrari, C. Landim, H. Thorisson}

\address{IME USP, Caixa Postal 66281, 05311-970 - S\~{a}o Paulo, Brasil.
  \newline \rm \texttt{pablo@ime.usp.br}, \href{http://www.ime.usp.br/~pablo}{http://www.ime.usp.br/\~{}pablo}}

\address{IMPA, Estrada Dona Castorina 110,
CEP 22460 Rio de Janeiro, Brasil and CNRS UMR 6085,
Universit\'e de Rouen, 76128 Mont Saint Aignan, France.
\newline
\rm \texttt{landim@impa.br},
\href{http://www.impa.br/Pesquisadores/Claudio/}{http://www.impa.br/Pesquisadores/Claudio/}
}

\address{ICE-UICE-SI,
         Science Institute,
         University of Iceland (H\'ask\'oli Islands),
         107 Reykjav\'{\i}k,
         Iceland\newline
\rm {\texttt hermann@hi.is},
\href{http://www.hi.is/~hermann/}{http://www.hi.is/\~{}hermann/}}

\maketitle

\paragraph{\bf Abstract}
We give a deterministic algorithm to construct a graph with no loops
(a tree or a forest) whose vertices are the points of a
$d$-dimensional stationary Poisson process $S\subset\R^d$. The
algorithm is independent of the origin of coordinates.  We show that
(1) the graph has one topological end ---that is, from any point there
is exactly one infinite self-avoiding path; (2) the graph has a unique
connected component if $d=2$ and $d=3$ (a tree) and it has infinitely
many components if $d\ge 4$ (a forest); (3) in $d=2$ and $d=3$ we
construct a bijection between the points of the Poisson process and
$\Z$ using the preorder-traversal algorithm.  To construct the graph
we interpret each point in $S$ as a space-time point
$(x,r)\in\R^{d-1}\times\R$. Then a $(d-1)$ dimensional random walk in
continuous time continuous space starts at site $x$ at time $r$. The
first jump of the walk is to point $x'$, at time $r'>r$, $(x',r')\in
S$, where $r'$ is the minimal time after $r$ such that $|x-x'|<1$.
All the walks jumping to $x'$ at time $r'$ coalesce with the one
starting at $(x',r')$. Calling $(x',r') = \alpha(x,r)$, the graph has
vertex set $S$ and edges $\{(s,\alpha(s)), s\in S\}$. This enables us
to shift the origin of $S^{\circ} = S + \delta_0$ (the Palm version of S) to
another point in such a way that the distribution of $S^{\circ}$ does not
change (to any point if $d = 2$ and $d = 3$; point-stationarity).



\section{Introduction}
Let $S^{\circ}$ be the Palm version of a stationary Poisson process $S$ in
$\bb R^d$, that is, $S^{\circ}$ is a random set with the same distribution
as $S \cup \{0\}$.  In the open problem session at the Brazilian
School of Probability in 2001 the third author posed the
following three problems. Let $\cS$ be the support of $S$ and $\cS^\circ$ the
support of $S^\circ$.

\begin{problem}
\label{p1}\rm
When $d > 1$, is there some non-randomized way of shifting the origin of
$S^{\circ}$ from the point at the origin to another point $X\in S^{\circ}$ so
that the distribution of $S^{\circ}$ as seen from this point remains the same?
More precisely, is there a map $\pi:\cS^\circ\to\R^d \setminus \{0\}$ such
that with $X = \pi(S^{\circ})\in S^{\circ}\setminus \{0\}$ and, with $S^{\circ} - X =
\{s - X : s \in S^{\circ}\}$,
\begin{equation}
  \label{(1)}
 S^{\circ} - X  \;=\;  S^{\circ} \quad\text{in distribution ?}
\end{equation}
This is clearly possible when $d = 1$, since then the intervals between points
are i.i.d. exponential and remain so when the origin is shifted to the $n$th
point on the right (or on the left) of the point at the origin.  When $d > 1$,
an $X$ satisfying \reff{(1)} ---and with $P(X \not = 0)$ arbitrarily close to
1--- exists if external randomization is allowed. This is shown by Thorisson
(1999 and Chapter 9 in 2000) who also proved that \reff{(1)} holds if the
\emph{point shift} map $\theta_{\pi}:\cS^\circ\to\cS^\circ$ defined by
\[
\theta_{\pi}(S^{\circ})\,:=\, S^{\circ} - \pi(S^{\circ})\,=\, S^{\circ} - X
\]
is \emph{bijective}.
\end{problem}

\begin{problem}
\label{p2}\rm
Does there exist a family of maps $\pi_n\,:\,\cS^{\circ}\to \R\setminus\{0\}$
such that, defining $X_n:=\pi_n(S^{\circ})$, $n\in \Z$, we have $X_n\neq
X_{n'}$ if $n\neq n'$ and
\begin{equation}
  \label{3}
  S^{\circ} = \{X_n : n \in \Z\}
\end{equation}
(that is, a labelling of all the points of $S^{\circ}$) and the maps
$\theta_{\pi_n}:\cS^\circ\to\cS^\circ$ defined by
\[
\theta_{\pi_n} (S^{\circ})\,:=\, S^{\circ}-X_n
\]
are bijective? As we just have seen, this would imply that $S^{\circ}-X_n =
S^{\circ}$, in distribution.
\end{problem}

\begin{problem}\label{p3}\rm
In the references mentioned before, Thorisson defined
``point-stationarity'' of a point process as distributional invariance
under bijective point-shifts ``against any independent stationary
background''; this concept is shown to be the characterizing property
of the Palm version of any stationary point process $S$ in $\R^d$. A
natural question is whether ``against any independent stationary
background'' can be removed from the definition. In other words, can
the definition of ``point-stationarity'' be reduced to distributional
invariance under non-randomized bijective point-shifts?
\end{problem}

Olle H{\"a}ggstr{\"o}m has given an example of a non-randomized bijection as
in Problem \ref{p1}: let $X$ be the closest point to the point at the origin
if the point at the origin is also the closest point to that point; otherwise
let $X$ be $0$.  But is there a strictly non-zero $X$? Dana Randall and the
first author proposed the following map in the same vein: in the first step
``marry'' each point $s$ to its closest point $s'$ if $s$ is the closest point
to $s'$. Call $M(S)$ the set of points in $S$ married using this procedure and
$M_1=M(S)$.  The set of points married in the $n$th step is $M_n =
M(S\setminus (M_1\cup\dots\cup M_{n-1}))$. Every point will eventually get
married, that is $\cup_n M_n = S$, because the Poisson process has no
descending chains, as proved by H{\"a}ggstr{\"o}m and Meester (1996); this has
been observed by Holroyd and Peres (2003). In this case, the map $s\mapsto$
(spouse of $s$) for $s=0$ is such a non-zero $X$ for all $d$.

We provide another solution to \ref{p1} for all $d$ and show that the answer
to \ref{p2} is {\bf ``yes''} when $d = 2$ and $d = 3$. We give a partial
response to \ref{p3}.

We show that in $d = 2$ and $d = 3$ it is possible to join the points of
$S^{\circ}$ in an origin-independent way into a connected tree with finite
branches. For this tree every pair of vertices have an ancestor in common and
every vertex has a finite number of descendents. The resulting tree is called
\emph{Poisson tree}. We order sets of sisters using the first spacial
coordinate and associate to each vertex a semi infinite sequence indicating
the sister-order of her ancestors. Then we order these sequences
lexicographically from the past to get a ``unique infinite succession line''
of vertices. The total order so obtained corresponds to the \emph{preorder
  traversal} algorithm used in computer science.

Let $\pi_0(S^{\circ}):=0$ and inductively $\pi_n(S^{\circ}):=X_n$, where $X_n$
is the successor of $X_{n-1}$ for $n\ge 1$ and the predecessor of $X_{n+1}$
for $n\le 1$ in the total order just described. Let $\theta_{\pi_n} :
S^0\mapsto S^0-X_n$. It is clear that each point is the successor of its
predecessor and more generally the $n$th successor of its $n$th predecessor
and vice versa:
\[ \pi_{-n}(S^{\circ}-\pi_n (S^{\circ}) ) = -\pi_n (S^{\circ} )\]
%
% \[
% \pi_{-n}(S^{\circ}-\pi_n(S^{\circ})) = \pi_n(S^{\circ}-\pi_{-n}(S^{\circ})) = 0,
% \]
for $n\ge 0$, which implies $\theta_{\pi_n}^{-1} = \theta_{\pi_{-n}}$ for all
$n$. Thus, for all $n\in\Z$, $\theta_{\pi_n}$ is a bijection as required in
\ref{p2}.  Notice also that $\theta_{\pi_n} = \theta_{\pi_{1}}^n$ the $n$th
iteration of $\theta_{\pi_1}$. In $d\ge 4$ we construct a ``forest'',
infinitely many connected trees with finite branches. Each tree can be ordered
using the preorder traversal algorithm and the map $s\mapsto$ (successor of
$s$) is a solution to Problem \ref{p1} in any $d\ge 1$.

In $d=2$ (one dimension for space and the other for time) our construction is
a continuous space-time analogous of a discrete space-time system of
coalescing random walks. Both the system of coalescing random walks and the
Poisson tree converge to the so called ``Brownian web'', a system of
one-dimensional coalescing Brownian motions starting at every space-time point
in $\R^2$. The random walk convergence and different properties of the web
have been studied by Arratia (1979, 1979a) and Toth and Werner (1998).
Ferrari, Fontes and Wu (2003) show that the Poisson tree converges to the
Brownian web in the sense proposed by Fontes, Isopi, Newman and Ravishankar
(2003).

Gangopadhyay, Roy and Sarkar (2002) proposed a system of coalescing random
walks based on a Bernoulli product measure in $\Z^2$ and show that it produces
a connected tree. Their motivation is to provide a model for drainage
networks; see the book of Rodriguez-Iturbe and Rinaldo (1997).

A tree having $S$ as vertex set is the \emph{minimal spanning tree}
constructed as follows. Choose arbitrarily an initial point $s\in S$.  Let
$(V_n,E_n)$ be the set of vertices and edges chosen up to the the $n$th
iteration; $(V_0,E_0)=(\{s\},\emptyset)$.  The $(n+1)$th point $s_{n+1}$ is
the point in the complement of $V_n$ that is closest to $V_n$. The $(n+1)$th
edge is the pair $(s',s_{n+1})$, $s'\in V_n$ realizing the distance between
$V_n$ and $s_{n+1}$. Alexander (1995) proved that in $d=2$ the construction is
independent of the initial point, that $V_n\to S$ and that the resulting tree
has all branches finite. In this terminology, the tree has one topological end
---from each point there is only one infinite self-avoiding path. It is
believed that for $d\le 8$ this tree has a unique connected component; Newman
and Stein (1992, 1994).

In Section \ref{s2} we relate trees and succession lines. In Section \ref{s3}
we construct the Poisson tree and state that in $d=2$ and $d=3$ the tree is
connected and has finite branches (Theorem \ref{20}). In Section \ref{s4} we
introduce a graphical construction of a system of coalescing random walks,
associate it to the tree and prove Theorem \ref{w20} ---from where Theorem
\ref{20} follows. In that section we also show that the system of coalescing
random walks is ergodic and converges exponentially fast to the unique
invariant measure. In Section \ref{s5} we give a limited reply to Problem
\ref{p3}. In Section \ref{s6} we make some final remarks and state related
open problems.

\section{Trees and succession lines}
\label{s2}

Let $\G$ be an oriented graph with no loops (a tree or a forest) such that
each vertex $s$ of $\G$ is the startpoint of exactly one outgoing edge and the
endpoint of either none or a finite number of ingoing edges.  The endpoint of
the outgoing edge is called the \emph{mother} of $s$, while the startpoints of
the ingoing edges are the \emph{daughters} of $s$. If two vertices are
daughters of the same vertex, they are \emph{sisters} of each other. We order
groups of sisters according with some spacial property like first coordinate,
distance to the mother, etc; we denote $s<s'$ and say that $s$ is older than
$s'$ in the adopted order. This order is not necessarily coherent with the
partial order induced by the ancestor-descendent relation.

Call $\alpha(s,\G)$ the mother of $s$.  Let $\alpha^0(s,\G) = s$, and
iteratively, for $n\ge 1$, $\alpha^{n}(s,\G) = \alpha(\alpha^{n-1}(s,\G),\G)$
the $n$th ancestor of $s$.  Let
\begin{eqnarray}
  \label{31}
  D^1(s,\G) &:=&\{s'\in S\,:\, \alpha(s',\G) = s\};\nonumber \\
  D^n(s,\G) &:=&  \{s'\in S\,:\,\alpha(s',\G) \in
  D^{n-1}(s,\G)\} \nonumber \\
D(s,\G) &:=& \bigcup_{n\ge 0} D^n(s,\G)
\end{eqnarray}
be respectively the first generation, the $n$th generation and the set
of all descendents of~$s$; we call $D(s,\G)$ the \emph{branch} of $s$
and say that $\G$ has finite branches if $D(s,\G)$ is finite for all
$s$. We say that two vertices $s$ and $s'$ are \emph{connected} if
they have an ancestor in common: there exist nonnegative integers $n$
and $m$ such that $\alpha^n(s,\G) = \alpha^m(s',\G)$. This defines an
equivalence relation in $S$; the equivalence classes are called
\emph{connected components}.

Let $\G$ be a connected tree with finite branches. Let $\sigma(s)=1$ if $s$ is
the eldest among her sisters, $\sigma(s)=2$ if $s$ is the second sister, and
so on. We associate to each vertex $s$ the sequence of relative sister-order
of its ancestors: let $\sigma_i(s):= \sigma(\alpha^i(s))$, $i\ge 0$.  If
$s''=\alpha^i(s) =\alpha^j(s')$ is the closest common ancestor of $s$ and
$s'$, for nonnegative $i,j$, then they can be lexicographically ordered using
the non common part of the sequences: we say that $s$ \emph{precedes} $s'$ if
$(\sigma_i(s),\sigma_{i-1}(s),\dots,\sigma_1(s),\sigma_0(s))$ is
lexicographically before
$(\sigma_j(s'),\sigma_{j-1}(s'),\dots,\sigma_1(s'),\sigma_0(s'))$. That is, if
$\sigma_{i-1}(s)< \sigma_{j-1}(s')$, with the convention $\sigma_{-1}(s) = 0$.

For $s\in\G$, define the \emph{successor} of $s$ as $s'\in\G$ if $s$ precedes
$s'$ and there is no vertex preceded by $s$ and preceding $s'$. Conversely $s$
is the \emph{predecessor} of $s'$ if and only if $s'$ is the successor of $s$.

The successor of a vertex can be found using the following algorithm. If the
vertex has a daughter, choose the eldest daughter. If it does not
have a daughter but has a younger sister, choose the eldest among its younger
sisters. If it does not have a daughter and not a younger sister, move up the
tree until you hit the first point that has a younger sister and choose the
eldest among its younger sisters. This requires that every vertex has an
ancestor with a younger sister.

The predecessor vertex can be found with this algorithm. If the vertex has an
elder sister, choose the youngest among her elder sisters and then move from
her down the tree choosing the youngest daughter in each step until you come
to a point with no daughter; this will be the predecessor. This requires that
the branch of the mother of each individual is finite.

We say that there is a \emph{succession line} from $s$ to $s'$ if there exists
a finite sequence of vertices $s = s_0,\dots, s_k = s'$ such that $s_{\ell-1}$
is successor of $s_{\ell}$ for $\ell=1,\dots,k$.  We say that the tree has an
\emph{infinite succession line} if every vertex has a predecessor and a
successor and that it has a \emph{unique infinite succession line} if
furthermore for every couple of vertices $s$, $s'$ there is a succession line
either from $s$ to $s'$ or from $s'$ to $s$. The following lemma follows from
the definitions.

\begin{lemma}
\label{w13}
If $\G=(S,E)$ has a unique infinite succession line, then the map $\pi(s,\G) =
s'$, the successor of $S$ and $\pi^{-1}(s',\G) = s$, the predecessor of $s'$
are well defined (i.e., both the successor and predecessor algorithms find a
vertex) and one is the inverse of the other. Furthermore, for all $s\in S$,
$S= \{\pi^n(s,\G): n\in \Z\}$.
\end{lemma}

The following lemma gives an equivalent condition.

\begin{lemma}
\label{w12}
A tree $\G$ has a unique connected component, finite branches and every vertex
has a mother and an ancestor with a younger sister if and only if $\G$ has a
unique infinite succession line.
\end{lemma}

\proof Since every vertex $s'$ has a mother and all branches are finite, the
predecessor of $s'$ is one of the vertices of the branch of her mother and can
be found in a finite number of steps. The condition that every vertex $s$ has an
ancestor with a younger sister guarantees the existence of the successor of
$s$ that can also be found in a finite number of steps. Since the tree is
connected, two arbitrary vertices $s$ and $s'$ have an ancestor in common, say
$s''$. The branch of $s''$ is finite, by hypothesis; say it has $n$ vertices.
Call $s''_0=s''$, $s''_k$ the successor of $s''_{k-1}$, $k=1,\dots,n$. Then
the branch of $s''$ is the same as $\{s''_0,\dots,s''_n\}$. In particular it
contains $s$ and $s'$ and a succession line either from $s$ to $s'$ or from
$s'$
to $s$. The converse statement is immediate. \square

\section{Poisson tree}
\label{s3}

Let $d\ge 2$ and $S$ be a locally finite configuration of points in $\R^d$.
For each $s=(s_1,\dots,s_d)\in \R^d$ call its first $d-1$ coordinates
$x(s)=(s_1,\dots,s_{d-1})\in\R^{d-1}$ and the remaining coordinate $r(s) =
s_d$. In this way $s=(x(s),r(s))$; $x(s)$ is interpreted as the \emph{space
coordinate} and $r(s)$ the \emph{time coordinate} of $s$.

For each $x\in\R^{d-1}$ let $B(x) = \{x'\in\R^{d-1};\, |x'-x|\le 1\}$ be the
$d-1$ dimensional Euclidean ball of radious 1 centered at $x$.  For each
$(x,r)\in S$ call $\{(x',r): x'\in B(x)\}$ (the $d-1$ dimensional disk
centered at $(x,r)$ perpendicular to the $d$th axis) the \emph{obstacle}
associated to $(x,r)$.  The set of obstacles is given by
\[{\cO}(S) = \bigcup_{(x,r)\in S} \{(x',r): x'\in B(x)\}\]

Think that each point $(x,r)\in S$ emits a laser ray in the positive $d$th
coordinate that is stopped by the obstacles. The first obstacle hit by the
ray of $(x,r)$ has second coordinate
\begin{equation}
   \label{w11}
   \tau((x,r),S):= \inf\{t>r\,:\,(x,t)\cap \cO(S)\neq\emptyset\}
\end{equation}
 with center
\begin{equation}
  \label{10}
  \alpha ((x,r),S)\;:=\; (x',r')\in S \quad\hbox{ if }\quad \tau(s,S) = r'
\end{equation}
which is called the \emph{mother of $(x,r)$}. Reciprocally, $s$ is a
\emph{daughter} of $\alpha(s,S)$. The above objects are well defined
for $S$ if $\tau(s,S)<\infty$ for all $s\in S$ and if no point has two
mothers.  In this case let $\G(S)=(S,E(S))$ be the random directed
graph with vertices $S$ and edges $E(S)=\{(s,\alpha(s,S))\,:\, s\in
S\}$. This graph has no loops, hence it is an oriented tree. Notice
that $\alpha(s,S)$ coincides with $\alpha(s,\G(S))$, in the notation
of the previous section.

\begin{thm}
  \label{20}
  Let $S$ be the realization of a homogeneous $d$-dimensional Poisson process
  and $S^{\circ}$ its Palm version.  Then, for $\G=\G(S)$ and $\G=\G(S^{\circ})$ it holds
  $S$-a.s. and $S^{\circ}$-a.s.:
  \begin{itemize}
  \item [(a)]$\G$ is well defined.
  \item[(b)] In $d=2,3$, $\G$ has a unique connected component.
  \item[(c)] In $d\ge 4$, $\G$ has infinitely many connected components.
  \item[(d)] All branches of $\G$ are finite.
  \item[(e)] Every vertex has a mother.
  \item[(f)] Using the order of the first coordinate in $\R^d$, each vertex
    has an ancestor with a younger sister.
  \end{itemize}
\end{thm}

Items (b) and (d) to (f) are necessary to construct the maps $\pi_n$
of Problem~\ref{p2}. To prove (a) it suffices to see that for almost
all realizations of $S$ and $S^{\circ}$, \reff{10} has a unique
solution; that is, each point has a unique mother. This is clear for a
Poisson process and its Palm version. The proofs of (b) to (e) are
based on a particle system in $\R^{d-1}$ studied in next section; they
are direct consequences of Theorem \ref{w20}.

\section{Coalescing random walks.}
\label{s4}
Let $S$ be a point configuration satisfying that every point has
exactly one mother. Write $\alpha(s)$ instead of $\alpha(s,S)$ and
recall the notation introduced in section 2.  For $s\in \R^d$, let
\begin{equation}
  \label{33}
  \tau^n(s) := r (\a^n(s))\,;
\end{equation}
that is, $\tau^n(s)$ is the time
coordinate of the $n$th ancestor of $s$.  For each $(x,r)\in \R^{d-1}\times\R$
 let $X^{(x,r)}_r=x$ and for $t\ge r$,
\begin{equation}
  \label{32}
  X^{(x,r)}_t = x(\a^n(x,r)),\qquad \hbox{ for }
t\in[\tau^n(x,r),\tau^{n+1}(x,r))\; ,\; n\ge 0.
\end{equation}
That is, $X^{(x,r)}_t$ starts at $x$ at time $r$ and remains still between
$\tau^{n-1}(x,r)$ and $\tau^{n}(x,r)$, when it jumps to $x(\alpha^{n}(s))$,
the center of the $n$th obstacle it meets. In other words, $\tau^n(x,r)$ is
the instant of the $n$th jump of the point that at time $r$ was at position
$x$. The family
\begin{equation}
  \label{w21}
\mc X(S):=  \{(X^{(x,r)}_t\,:\,t\ge r): (x,r)\in S\}
\end{equation}
is a deterministic function of $S$.

Let $S$ be a point configuration of a homogeneous Poisson process of rate
$\lambda$. Then \reff{w21} is a family of random processes in the probability
space where $S$ is defined; its law corresponds to a system of coalescing
random walks with births, whose marginal distribution is described by ``each
random walk waits an exponentially distributed random time of mean $(V_{d-1}
\la)^{-1}$ after which it chooses a point uniformly in $B(x)$, the
$(d-1)$-dimensional Euclidean ball of radius one centered at $x$, and jumps to
it''. Here $V_{d}$ stands for the volume of the $d$-dimensional ball of radius
one.  Particles are created at a Poisson rate $\lambda$ and the interaction
appears when two walks are located at points $x$ and $y$, $|x-y|<2$ and a
Poisson event appears in $B(x)\cap B(y)$; in this case both walks jump to the
same point and coalesce.  Disregarding the label of the particles, for $t\in
\R$, let
\[
\eta_t=   \{X^{(x,r)}_t\,:\, (x,r)\in S, r\le t\}\; .
\]
Here $\eta_t$ is a discrete subset of $\R^{d-1}$. The process starting from a
fixed configuration $\eta\subset\R^{d-1}$ at fixed time $t'$ is defined by:
\[
\eta^{\eta,t'}_t = \{X^{(x,r)}_t\,:\, (x,r)\in S,\, t'\le r\le t \, \text{ or
} x\in \eta,\, r=t'\}\;
\]
for $t\ge t'$. That is, ignore the starting points with second coordinate less
than $t'$ and include the points with first coordinate in $\eta$ and second
coordinate equal to $t'$.

\begin{propos}
  \label{39}
  Let $S$ be a stationary Poisson process.  The process
  $(\eta^{\eta,t'}_t\,:\,t\ge t')$ is Markov with generator defined by
\begin{equation*}
  Lf(\eta) = \la \int_{\R^{d-1}}
  [f(\eta\setminus B(u)\cup\{u\})-f(\eta)] du\,,
\end{equation*}
for functions $f:\cS_{d-1}\to\R$ depending on bounded regions, where
$\cS_{d-1}$ is the set of locally finite labeled configurations of $\bb
R^{d-1}$.  Furthermore $\eta_t$ is a stationary version of the unique
stationary process with generator $L$. In particular, the marginal
distribution of $\eta_t$ for any given $t$ is the unique invariant
measure for the generator $L$.
\end{propos}

\proof The independence properties of the Poisson process $S$ imply that for
any initial configuration $\eta$, the process $(\eta^{\eta,t'}_t\,,\;t\ge t')$
is Markov and has generator $L$.

To show that the process is well defined starting at time $-\infty$ we
show first that the configuration in finite sets depends only on a
finite number of points of $S$.  Let $\Lambda$ be a subset of $\R^{d-1}$
with finite Lebesgue measure. The idea is to prove that for almost all
configuration $S$, $\eta_t\cap \Lambda$ depends only on a finite (but
random) subset $S_{\Lambda,t}$ of $S$.  To get the set $S_{\Lambda,t}$ one
first
translates the set $\Lambda$ backwards in time hitting points $(x,r)$. Each
time a point $(x,r)$ is hit, $\Lambda$ is updated to $\Lambda\setminus
B(x)$ and
$(x,r)$ is included in $S_{\Lambda,t}$. The procedure stops when $\Lambda$ is
updated to the empty set. Then for each point $(x,r)$ already in
$S_{\Lambda,t}$ include also all its ancestors with second coordinate less
than $t$: $(x^i,r^i)= \alpha^i(x,r)$ with $r^i\le t$. We leave to the
reader to show that $\eta_t\cap \Lambda$ depends only on $S_{\Lambda,t}$.

For each finite measure set $\Lambda$, $S_{\Lambda,t}$ has a finite number of
points with probability one. This shows that for any fixed $t$, as
$t'\to-\infty$ the variable $\eta^{\eta,t'}_t\cap \Lambda$ converges almost
surely to a random variable $\eta_t\cap \Lambda$. Since the law of
$\eta_t\cap \Lambda$ does not depend on $t$, it is invariant for the
process.  \square

It follows from this construction that the process converges
exponentially fast to equilibrium:

\begin{corol}
  Denote by $\{P_t : t\ge 0\}$ the semigroup associated to the
  generator $L$. For every bounded function $f: \cS_{d-1} \to \R$
  with finite support, there exists a finite constant $C(f)$ such that
$$
\big\Vert P_t f - E[f] \big\Vert_\infty \;\le\; C(f) e^{- C_d\lambda
  t}
$$
for every $t\ge 0$ and some finite constant $C_d$ depending only on
the dimension. In this formula $E[\, \cdot\,]$ stands for the
expectation with respect to the stationary state.
\end{corol}

\proof Fix two configurations $\eta$, $\xi$ and a finite cube $\Lambda$. For
$t\ge 0$, denote by $a_t^\Lambda$ all points $x$ in $\Lambda$ such that $(x,r)$
belongs to $S$ for some $0\le r\le t$ and denote by $A_t^\Lambda$ the union
of all $(d-1)$-dimensional balls of radius one with center in $a_t^\Lambda$:
\begin{eqnarray*}
&& a_t^\Lambda\;=\; \bigcup \big\{ x : (x,r)\in S \quad\text{for some $0
\le r\le t$ }\big\} \;, \\
&&\quad A_t^\Lambda \;=\; \bigcup_{x\in a_t^\Lambda} B(x) \; .
\end{eqnarray*}
For a configuration $\zeta$, denote by $\eta^\zeta_t$ the state at
time $t$ of the process which starts from $\zeta$.  By construction
$\eta_t^\eta \cap \Lambda = \eta_t^\xi \cap \Lambda$ if $A_t^\Lambda$
contains $\Lambda$.

Divide the cube $\Lambda$ in small cubes, in such a way that if each small
cube contains a point in $a_t^\Lambda$ then $A_t^\Lambda$ contains
$\Lambda$.  Denote
these cubes by $\{E_i : 1\le i \le M\}$ and notice that the number $M$
of them is equal to $C |\Lambda|$ for some constant $C$ depending only on
the dimension. In view of the two previous observations,
$$
\P[ \eta_t^\eta \cap \Lambda \not = \eta_t^\xi \cap \Lambda] \;\le\;
\P[\Lambda \not \subset A_t^\Lambda] \;\le\; \P\Big[ \bigcup_{i=1}^M
\Omega_i(t) \Big]\; ,
$$
where $\Omega_i(t)$ is the event that the cube $E_i$ does not contain
points in $a_t^\Lambda$. Since a point appears in a set $A$ at rate $\lambda
|A|$, $P[\Omega_i(t)] \le \exp\{-\lambda |E_i| t\}$ so that
$$
\P[ \eta_t^\eta \cap \Lambda \not = \eta_t^\xi \cap \Lambda] \;\le\;
C_1 |\Lambda| e^{- C_2 \lambda t}
$$
for two finite constants depending only on the dimension because $M =
C(d) |\Lambda|$ and $|E_i| = C(d)$.

To conclude the proof, it remains to consider a bounded function $f$
with finite support. Denote by $\Lambda$ a cube which contains its support.
$\Vert \, P_t (f) - E[f]\, \Vert_\infty$ is bounded by
$$
\sup_{\eta,\xi} \E \Big[ \, \big \vert f(\eta^\eta_t) - \E
[f(\eta^\xi_t)] \big \vert\, \Big ] \;\le\; 2 \, \Vert f\Vert_\infty
\, \P[\eta_t^\eta \cap \Lambda \not = \eta_t^\xi \cap \Lambda] \; ,
$$
which concludes the proof of the lemma. \square

\begin{thm}
  \label{w20}
  Let $S$ be the realization of a homogeneous $d$-dimensional Poisson process
  and $S^{\circ}$ its Palm version.  Then, for $\mc X=\mc X(S)$ and $\mc X=\mc
  X(S^{\circ})$ it holds $S$-a.s. and $S^{\circ}$-a.s.:
  \begin{itemize}
  \item[(a)] $\mc X$ is well defined.
  \item[(b)] In $d=2,3$, every couple of walkers $X^{x,r}_t$, $X^{x',r'}_t$ in
    $\mc X$ will eventually meet.
  \item[(c)] In $d\ge 4$, there are infinitely many walkers in $\mc X$ that do
    not meet.
  \item[(d)] Each walker alive at time $t$ was born at a finite time
    before $t$.
  \item[(e)] Every walk will eventually jump.
  \item[(f)] Every walk will eventually coalesce with a younger walk (in
    the order of the first coordinate at coalescence time).
  \end{itemize}
\end{thm}

\paragraph{\bf Proof.} (a) and (e) follow from the properties of the Poisson
process: Each walk will eventually hit a unique obstacle and jump to it.

By symmetry each walk $X^{(x,r)}_{t}$ is a martingale: Letting $\cF_t$
be the $\sigma$-algebra generated by $\{s\in S\,:\, s_d\le t\}$,
\begin{equation*}
  \E(X^{(x,r)}_{t'}\,|\, \cF_t) \,=\, X^{(x,r)}_t,\qquad t'\ge t\ge
  r\; .
\end{equation*}

\vskip 3mm 

(b) In $d=2$, for $x'<x$, the difference walk $D_t:=X^{(x,r)}_{t}-
X^{(x',r')}_{t}\in\R$, $t\ge\max\{r,r'\}$ is a positive martingale. To
check it is enough to show that $\partial_t \E(D_{t}
\,|\,X^{(x,r)}_{t} =a, \, X^{(x',r')}_{t}=b) =0$, which follows from
an elementary computation divided in two cases: $a<b-1<a+1<b$ and
$a+1<b-1$.  Hence $D_t$ it is recurrent in the sense that $D_t<1/2$
will occur for an infinite number of times $t$. Since each time
$D_t<1$ there is a positive probability that $D_{t+1}=0$, and 0 is an
absorbing point for $D_t$, this implies that $(x,r)$ and $(x',r')$
have an ancestor in common. Since this follows for all points, there
is only one tree in $d=2$.

For $d=3$, let $Y_t$ and $Y'_t$ be two $(d-1)$-dimensional independent
random walks with the same marginal distribution as
$X_t=X^{(x,r)}_{t}$ and $X'_t=X^{(x',r')}_{t}$ respectively. Without
loss of generality we can assume $t\ge 0$. The generator $\mc L$ of
the process $(Y_t, Y'_t)$ is given by
$$
(\mc L f)(y,y') \; =\; \la \int_{B(y)} [f(z,y') - f(y,y')] \, dz
\; +\; \la \int_{B(y')} [f(y,z) - f(y,y')] \, dz \; .
$$
Since in dimension $2$, $\log |z|$ is a harmonic function, $\mc L
\log |y-y'| = 0$ so that $\log |Y_t- Y'_t|$ is a local martingale.
This observation and standard arguments involving the hitting time of
two balls centered at the origin and of radius $\uR < |y-y'| < \oR$
permits to show that $Y_t- Y'_t$ is recurrent.

Assume $|X_0-X'_0|>2$ and couple $((X_t,X'_t),(Y_t,Y'_t))$ up to the
first time $T$ such that $|Y_t- Y'_t|<2$; for $t\in [0,T]$ we have
$X_t=Y_t$ and $X'_t=Y'_t$. Then wait up to the first time $T'>T$ such
that either $X_{T'}=X'_{T'}$ or $|X_{T'}-X'_{T'}|>2$. Since there is a
positive probability that $X_{T'}=X'_{T'}$ and the recurrence of
$Y_t- Y'_t$ guarantees the existence of infinitely many attempts,
$X_{t}=X'_{t}$ eventually with probability one.

\vskip 3mm

(c) For $d\ge 4$, consider the process $(Y_t, Y'_t)$ defined in (b).
By similar reasons to the ones presented in (b), $|Y_t- Y'_t|^{2-d}$
is a local martingale. This observation and standard arguments
involving the hitting time of two balls centered at the origin and of
radius $\uR < |y-y'| < \oR$ permits to show that the process $Y_t-
Y'_t$ is transient.

Fix $\uR>2$ and proceed as in the proof of (a) and couple the
independent and the interacting walks in such a way that they coincide
if the distance is bigger than two. Hence, for initial points $y$,
$y'$ with $|y-y'|>2$ we get
\begin{eqnarray}
  \label{w15}
\P(|X^{y,0}_t - X^{y',0}_t| =0\text{ for some }t\ge 0)
&\le &
\P(|X^{y,0}_t - X^{y',0}_t|  <2 \text{ for some }t\ge 0)\nonumber\\
&\le &
\P(|Y^{y,0}_t - Y^{y',0}_t|  <2 \text{ for some }t\ge 0)\nonumber\\
 &\le & \frac {\uR^{d-2}}{|y-y'|^{d-2}}\cdot\nonumber
\end{eqnarray}
Last inequality is obtained from the fact that $|Y_t- Y'_t|^{2-d}$ is
a local martingale and standard estimates involving hitting times of
balls centered at the origin.  Fix $\uR>2$, $n\ge 2$, $\epsilon>0$ and
$g(n,d,\epsilon,\uR)$, to be defined later.  Since for the Poisson
process we can always choose a point exterior to any bounded region,
choose $x_1, \dots, x_n\in \eta_0$ iteratively from $\eta_0$, a
configuration picked from the invariant measure for the coalescing
random walks, in the following way: Pick $x_1$ arbitrary, $x_{i}$ at
distance bigger than $g$ of the $i-1$ points chosen before,
$i=1,\dots,n$. Then
\begin{eqnarray}
  \label{w16}
 \lefteqn{ \P(\text{there are less than }n\text{ walkers that do not meet})}
\nonumber\\
&\le & \P(\exists i,j\in\{1,\dots,n\} \text{ with }|X^{x_i,0}_t -
X^{x_j,0}_t| =0\text{ for some }t\ge 0)\nonumber\\
&\le & \sum_{i,j}\P(|X^{x_i,0}_t - X^{x_j,0}_t| =0\text{ for some }t\ge
0)\nonumber\\
&\le & \frac{n(n-1)}{2}\frac {\uR^{d-2}}{|g(n,d,\epsilon,\uR)|^{d-2}}\;\le\;
\epsilon \nonumber
\end{eqnarray}
for a sufficiently large $g(n,d,\epsilon,\uR)$.  Hence
\[\P(\text{there are a finite number of  walkers that do not meet}) \;=\;0.\]

\vskip 3mm

(d) Let $\xi^r_t\subset \eta_t$ be the set of points of $\eta_t$ with branches
going up to time $r<t$. Obviously $\xi^r_t$ is a decreasing set in $r$:
$\xi^r_t\subset \xi^{r'}_t$ for $r\le r'\le t$.  Hence the set
$\xi_t=\lim_{r\to-\infty}\xi^r_t$ is well defined and describes the set of
points having infinite branches backwards in time. By construction the law of
$\xi_t$ does not depend on $t$ and it is stochastically bounded by the law of
$\eta_t$. The process $\xi_t$ is stationary by construction and Markovian
with generator
\begin{equation*}
  L_of(\xi) = \la \int_{\R^{d-1}}
  [f(\{\xi\setminus B(u)\}\cup\{u\})-f(\xi)]\, \one\{\xi\cap B(u)\neq
  \emptyset\}\, du\,,
\end{equation*}
for functions $f:\cS_{d-1}\to\R$ depending on bounded regions. It describes the
motion of coalescing random walks (without births).

We want to show that $\xi_t=\emptyset$ with probability one: $\P(\xi_t
=\emptyset)=1$ for all $t$ in $\bb R$.

Assume by contradiction that $\P(|\xi_t|\ge 1) >0$, where $|\xi|$
stands for the total number of particles in $\xi$.  We claim that
$\P(|\xi_t|=1)+ \P(\xi_t =\emptyset)<1$.  Indeed, if $\P(|\xi_t|=1)+
\P(\xi_t =\emptyset)=1$, for any finite set $A$, and any $t\in \R$,
$t'\ge 0$,
\begin{eqnarray*}
\P(\xi_{t+t'} \cap A = \emptyset) \;=\; \P(\xi_{t} = \emptyset) \;+\;
\P(\xi_{t+t'} \cap A = \emptyset \, |\, |\xi_t|=1 ) \P(|\xi_t|=1)\;.
\end{eqnarray*}
If $|\xi_t|=1$, for $t'\ge 0$, $\xi_{t+t'}$ is a symmetric random walk
on $\R^{d-1}$.  In particular, $\P(\xi_{t+t'} \cap A = \emptyset \, |\,
|\xi_t|=1 )$ converges to $1$ as $t'\uparrow\infty$. Hence, by
stationarity,
\begin{eqnarray*}
\P(\xi_{t} \cap A = \emptyset) \;=\; \P(\xi_{t} = \emptyset) \;+\;
\P(|\xi_t|=1) \; =\; 1
\end{eqnarray*}
for all sets $A$. Letting $A\uparrow \R^{d-1}$, we obtain that
$\P(\xi_{t} = \emptyset)=1$ in contradiction with the first
assumption.

We just proved that $\P(|\xi_t|=1)+ \P(\xi_t =\emptyset)<1$ so that
$\P(|\xi_t|\ge 2) >0$. In particular, there exists a bounded set $A$
such that $\P(|A\cap \xi_t| \ge 2)>0$.

Let $\Lambda$ be a $(d-1)$ dimensional cubic box centered at the origin. Let
$\beta_t= \E(|\xi_t\cap \Lambda|)$.  By invariance of the construction of
$\xi_t$, $\beta_t$ does not depend on $t$.  On the other hand, $\beta_t$ grows
at most at rate
\[
\la|\delta \Lambda| = \lambda O(|\Lambda|^{(d-2)/(d-1)})
\]
corresponding to the entrance of infinite branches through the boundary of
$\Lambda$. This happens at rate $\lambda$ times the volume of the set of
points in $\Lambda$ being at distance less than one from some point outside
$\Lambda$; this volume is of the order $|\Lambda|^{(d-2)/(d-1)}$.

On the other hand, $\beta_t$ decreases at least at rate
\[
\lambda \int_\Lambda \P(B(u)\cap \xi_t \ge 2)\, du \; =\;
\lambda \P(B(0)\cap \xi_t \ge 2) |\Lambda| \,.
\]
where we recall $B(u)$ is the $(d-1)$ euclidean ball of radius 1 centered at
$u$. The identity follows from the translation invariance of
$\xi_t$. We showed above that there exists a bounded set $A$ such that
$\P(|A\cap \xi_t| \ge 2)>0$. Since in a finite amount of time one can
find Poisson points $(x,r)$ taking the (at least) two points in $A$ at
distance less than one, $\P(B(u)\cap \xi_t \ge 2)>0$. Therefore,
$\beta_t$ decreases at least by a constant times $\lambda |\Lambda|$.
We conclude that
\begin{equation*}
  \frac{d}{dt}\beta_t =  O(|\Lambda|^{(d-2)/(d-1)})- O(|\Lambda|) <0
\end{equation*}
for sufficiently large $\Lambda$. This is in contradiction with the
time independence of $\beta_t$ and shows that $\P(\xi_t\neq \emptyset) =0$.

\vskip 3mm (f) For almost every point configuration $S$, two different points
in $S$ have all coordinates different. Hence they are distinctly ordered by
the first coordinate. We want to show that every walk will eventually coalesce
with a younger walk at coalescence time. Recall that $\tau_n(x,r)$ are the
jump times of the walk $X^{x,r}_t$. Let $K_n$ be the intersection of the
following events:
\begin{itemize}
\item Between $\tau_n$ and $\tau_n+1$, $X^{x,r}_t$ does not jump
\item Between $\tau_n$ and $\tau_n+1$ a new walk is born at distance bigger
  than one and less than two of $X^{x,r}_t$, call it $X^{x',r'}_t$;
  $|x'-X^{x,r'}_t|\in(1,2)$, $r'\in [\tau_n,\tau_n+1)$
\item Between $\tau_n+1$ and $\tau_n+2$ both walks coalesce.
\end{itemize}
$\tau_n$ is a hitting time: the event $\{\tau_n\le t\}$ is $\mc F_t$
measurable, i.e. it depends only on the points of the past of $t$. The above
event depends on points in a finite region in the future of $\tau_n$. It has
positive probability and given that it has occurred, the probability that
$x'\le X^{x,r'}_t$ and $x'\ge X^{x,r'}_t$ are equal. In the second case
$X^{x,r}_t$ has coalesced with a younger walk. Since $\tau_n$ is a Poisson
process, we can make infinitely many independent attempts, each with positive
probability of success. Hence with probability one each walk will coalesce
with a younger walk. \square

Notice that our proof of (f) uses the independence of disjoint regions of the
Poisson process and the fact that a success is attained after a geometric
number of independent attempts. Another possibility would be to use ergodic
theory \emph{a la} Burton and Keane (1989); see also Holroyd and Peres (2003).

\section{On point-stationarity}
\label{s5}
In this section we shall briefly consider Problem \ref{p3} from the
introduction.

Let $S^{\circ}$ be a point process in $\R^d$ with a point at the origin,
$0 \in S^{\circ}$.  According to Thorisson (1999, see also Chapter 9 in
2000), $S^{\circ}$ is \emph{point-stationary} if for any shift-measurable
stationary random field $Y = (Y_u : u \in \R^d)$, which is
independent of $S^{\circ}$, it holds that
\[
S^{\circ} - \pi(S^{\circ}, Y)  =  S^{\circ} \quad\text{in distribution,}
\]
where $\pi$ is such that the associated point-shift is a bijection.
It is further shown that
$S^{\circ}$ is point-stationary if and only if, for all $u \in \R^d$,
\begin{equation}
  \label{(a)}
 E\Bigl[f(u + S^{\circ} - U) \, |C^{\circ}| \Bigr] =
E\Bigl[f(S^{\circ} - U) \, |C^{\circ}| \Bigr] \; ,
\end{equation}
where $C^{\circ}$ is the Voronoi cell of the point at the origin and,
conditionally on $S^{\circ}$, $U$ is uniformly distributed on $C^{\circ}$.  In
particular when $E[\, |C^{\circ}|\, ] < \infty$, this means that the
reversed Palm version $S$ of $S^{\circ}$ is stationary,
\[
S - u = S \; ,\quad\text{in distribution for $u \in \R^d$.}
\]
We shall now show that if there exist non-randomized bijective point
shifts like those constructed in the introduction for the Poisson
process in the cases $d$ = 2 and 3, then point-stationarity
reduces to invariance under those point shifts. That is, the
stationary independent background field $Y$ is not needed to ensure
point-stationarity in those cases.

\begin{thm}
\label{thm}
Let $S^{\circ}$ be a point process in $\R^d$. Let
$\pi_n$, $n \in \Z$ be point maps and $\theta_{\pi_n}$
the associated point shifts.
Suppose that
\begin{equation}
  \label{(b)}
 \pi_{-n}\theta_{\pi_n} = -\pi_n\; ,
\end{equation}
that the points $X_n = \pi_n(S^{\circ})$, $n \in \Z$, are all distinct,
and that $S^{\circ} = \{X_n : n \in \Z\}$ a.s.  If
\begin{equation}
  \label{(c)}
S^{\circ} - X_n  =  S^{\circ}\; , \quad\text{in distribution for $n \in Z$,}
\end{equation}
then $S^{\circ}$ is point-stationary.
\end{thm}

\proof Let $C^n$ be the Voronoi cell of $X_n$. Due to \reff{(b)} and
\reff{(c)},
\[
\int_{(u + C^{\circ})\cap C_n}f(s + S^{\circ})\, ds = \int_{(u + C^{-n})\cap C_0}
f(s + S^{\circ})\, ds
\]
in distribution.  Take expectations and sum over $n$ to obtain
\[
\E\Bigl[\int_{u + C^{\circ}} f(s + S^{\circ})ds\Bigr]
\;=\; \E\Bigl[\int_{C^{\circ}}f(s + S^{\circ})ds\Bigr]\; .
\]
This is a reformulation of \reff{(a)} which is equivalent to
point-stationarity, and the proof is complete. \square

\section{Final remarks and open problems}
\label{s6}
We have deterministically constructed connected trees with a unique infinite
succession line having as vertices the points of a homogeneous Poisson process
in dimensions $d=2$ and $d=3$. This naturally poses the following problems.

(a) Construct an infinite succession line for a Poisson process in dimension
$d\ge 4$. Due to Theorem \ref{thm}, this would in particular solve Problem
\ref{p2}.

(b) Construct an infinite succession line independent of
the choice of a direction. The minimal spanning tree construction has this
property, but it is still to be proven that it possesses a unique infinite
succession line. On the other hand the minimal spanning tree has been proven
to be connected only in dimension 2.

Holroyd an Peres (2003) have solved (a) and (b) in all dimensions by
constructing a connected tree with finite branches having as vertices the
points of the Poisson process in a translation and rotationally invariant way.

A challenging problem is to give sufficient conditions for an ergodic
translation invariant point process in any dimension to be ordered in
a deterministic way in a unique infinite succession line. This would be an
important step towards solving Problem \ref{p3}.

The infinite succession line in $d=2$ corresponds to the ``random maze''
introduced in the Figure 2 of Toth and Werner (1998). Propositions 3.1 and
Lemmas 3.2 and 3.4 of that paper together with the convergence of the tree to
the Brownian web imply that the infinite succession line
converges to a random line that totally order the two-dimensional real
numbers.

\section*{Acknowledgements}
We thank Dana Randall and Xian-Yuan Wu for calling our attention to the
minimal spanning tree, Ken Alexander for his comments and references regarding
the minimal spanning tree, Balint T\'oth for enlighting discussions about
succession lines, Anish Sarkar for the references on drainage networks and
Yuval Peres for indicating the reference to H\"aggstr\"om and Meester.

This paper is partially supported by FAPESP, CNPq, PRONEX.

\end{document}